\documentclass{article}
\usepackage{amssymb}
\usepackage{amsmath}
\usepackage{cite}

\newtheorem{theorem}{Theorem}

\newtheorem{corollary}[theorem]{Corollary}

\begin{document}

\title{\ Lie and Noether point symmetries of a class of quasilinear systems
of second-order differential equations}
\author{Andronikos Paliathanasis\thanks{%
Email: anpaliat@phys.uoa.gr} \\
%EndAName
{\ \textit{Instituto de Ciencias F\'{\i}sicas y Matem\'{a}ticas, Universidad
Austral de Chile, Valdivia, Chile}} \and Michael Tsamparlis\thanks{%
Email: mtsampa@phys.uoa.gr} \\
%EndAName
{\ \textit{Faculty of Physics, Department of
Astronomy-Astrophysics-Mechanics,}}\\
{\ \textit{University of Athens, Panepistemiopolis, Athens 157 83, Greece}}}
\maketitle

\begin{abstract}
We study the Lie and Noether point symmetries of a class of systems of
second-order differential equations with $n$ independent and $m$ dependent
variables ($n\times m$ systems). We solve the symmetry conditions in a
geometric way and determine the general form of the symmetry vector and of
the Noetherian conservation laws. We prove that the point symmetries are
generated by the collineations of two (pseudo)metrics, which are defined in
the spaces of independent and dependent variables. We demonstrate the
general results in two special cases (a) a system of $m$ coupled Laplace
equations and (b) the Klein-Gordon equation of a particle in the context of
Generalized Uncertainty Principle. In the second case we determine the
complete invariant group of point transformations, and we apply the Lie
invariants in order to find invariant solutions of the wave function for a
spin-$0$ particle in the two dimensional hyperbolic space.\newline
\end{abstract}

\section{Introduction}

Lie symmetries is a powerful tool for the study of differential equations,
because they provide invariant functions which can be used to reduce the
order of a differential equation or reduce the number of variables, and
possibly lead to the determination of analytic solutions.\ For differential
equations which arise from a variational principle, i.e. follow from a
Lagrangian, the Lie point symmetries which in addition leave the action
invariant are called Noether point symmetries. Lie point symmetries span a
Lie algebra and their specialization Noether point symmetries span a
subalgebra. According to Noether's theorem to each Noether point symmetry
there corresponds a conservation law \cite{Stephani,Bluman,Katzin}.
Conservation laws play an important role in Classical Mechanics, General
Relativity, field theory and in the study of dynamical systems in general
\cite{TopSol,ExSol,Leach87,Leach88,Ibrag98,ibra2,Azad}.

In this work we study the Lie point symmetries and the Noetherian
conservation laws of the class of second-order differential equations which
follow from the Lagrangian\footnote{%
The Latin indices $i,j,...$ take the values $1,2,...,n,$ and the capital
indices $A,B,...$ take the values $1,2,...m.$\newline
}

\begin{equation}
L\left( x^{k},u^{C},u_{,k}^{C}\right) =\frac{1}{2}\sqrt{g}%
g^{ij}H_{AB}u_{,i}^{A}u_{,j}^{B}-\sqrt{g}V\left( x^{k},u^{C}\right)
\label{pd.01}
\end{equation}%
where $g_{ij}=g_{ij}\left( x^{k}\right) ,~g^{ij}g_{ij}=\delta _{j}^{i}$; $%
H_{AB}=H_{AB}\left( u^{C}\right) $, with $H_{AB}H^{AB}=\delta _{B}^{A}$, and
$u_{,i}^{A}=\frac{\partial u^{A}}{\partial x^{i}}$. $\ H_{AB}$ is the metric
of the $m$ dependent variables $u^{C}\left( x^{k}\right) $, $\dim H_{AB}=m$,
$g_{ij}$ is the metric of the $n$ independent variables~$x^{i}$,~~$\dim
g_{ij}=n$,~and $\delta _{\beta }^{\alpha }$ is the Kronecker delta.

Lagrangian (\ref{pd.01}) leads to the following system of Euler-Lagrange
equations%
\begin{equation}
P^{A}\left( x^{k},u^{C},u_{,i}^{C},u_{,ij}^{C}\right) \equiv
g^{ij}u_{,ij}^{A}+g^{ij}C_{BC}^{A}u_{,i}^{B}u_{,j}^{C}-\Gamma
^{i}u_{,i}^{A}+F^{A}\left( x^{k},u^{C}\right) =0,  \label{pd.02}
\end{equation}%
where $C_{BC}^{A}=C_{BC}^{A}\left( u^{D}\right) ,~\Gamma ^{i}=g^{jk}\Gamma
_{jk}^{i}\left( x^{r}\right) $ are the connection coefficients of the
metrics $H_{AB}$ and $g_{ij}$ respectively, and $F^{A}=H^{AB}V_{,B}$. The
system (\ref{pd.02}) consists of $m$ equations and depends on $n$ variables;
we call it a $n\times m$ system. For $n=1$ the system (\ref{pd.02}) reduces
to $m$ second-order ordinary differential equations, and for $m=1$ the
system (\ref{pd.02}) describes a second-order partial differential equation.

The determination of Lie point symmetries of the system (\ref{pd.02})
consists of two steps: (a) the determination of the conditions which the
symmetry vector must satisfy (symmetry conditions), and (b) the solution of
these conditions. The first step is formal, however the symmetry conditions
which arise can be quite involved. One way to \textquotedblleft
solve\textquotedblright\ the system of the symmetry conditions is to write
them in geometric form and then use the methods of Differential Geometry to
solve them. In this way the determination of the Lie point symmetries of a
differential equation is reduced to a problem of Differential Geometry where
there is an abundance of known results and methods to work. Indeed the Lie
point symmetries of the $1\times m$ systems of the form (\ref{pd.02}) have
been solved in this manner. Specifically, it has been proved that the Lie
point symmetries of a $1\times m$ system are generated by the elements of
the special projective algebra of the space $H_{AB}$ of the dependent
variables \cite{GNgrg,JPA2d}, and the Lie point symmetries form the
projective group of an affine space $V^{1+m}$ \cite{Aminova,Aminova2}.
Moreover for a class of singular $1\times m$ systems in which the
Hamiltonian function is vanished the Lie and the Noether point symmetries
follow from the conformal algebra of the space of the dependent variables
\cite{chris}.

The geometric approach has also been applied to the case of $n\times 1$
equations of the form (\ref{pd.02}) and it has been proved that in this case
the Lie point symmetries are generated by the elements of the conformal
algebra of the space of the independent variables $g_{ij}$ \cite%
{JGPpde,AnIJGMMP,Bozhkov}.

In the following we generalize the above results in the case of the $n\times
m$ system (\ref{pd.02}). In particular, we show that the Lie point symmetry
vectors in the space $\left\{ x^{i},u^{C}\right\} $ follow from the affine
collineations of the metric $H_{AB}$ and the conformal Killing vectors of
the metric $g_{ij}$. Moreover, there exists a connection between the two
algebras if and only if the space $H_{AB}$ admits a gradient homothetic
vector. The structure of the paper is as follows.

In section \ref{Preliminaries}, we present the basic definitions concerning
the Lie and the Noether point symmetries of differential equations as well
as the collineations of a Riemannian space. The Lie point symmetries of the
system (\ref{pd.02}) are studied in section \ref{Lie}, where we prove that
the generic Lie point symmetry vector is generated by the affine algebra and
the conformal algebra of the two metrics $H_{AB}$ and $g_{ij}$ respectively.
For the Noether point symmetries of Lagrangian (\ref{pd.01}) we derive the
generic form of the Noether vector and of the corresponding conservation
law. In section \ref{Laplace}, we apply the general results of the previous
sections to study the case of $m$ coupled Laplace equations. We find that if
the system of $m$ coupled Laplace equations (with $n>2$) admits $\frac{%
\left( n+2\right) \left( n+1\right) }{2}+m\left( m+1\right) $ Lie point
symmetries, then the two spaces $H_{AB}$,$~g_{ij}$ are flat.

In section \ref{GUPkg}, we consider the modified Klein-Gordon equation of a
particle in the Generalized Uncertainty Principle which is a fourth order
partial differential equation. With the use of a Lagrange multiplier we
reduce this equation to a system of two second-order partial differential
equations of the form of system (\ref{pd.02}). We study the Lie and the
Noether point symmetries of the new system and show that if the $g_{ij}$
space admits an $N$ dimensional Killing algebra then, the Lie and the
Noether point symmetries form Lie algebras of dimension $N+3$ and $N+2$
respectively. We apply this result in two cases of special interest: (A) the
particle lives in the flat Minkowski space-time $M^{4}$, and (B) the
particle lives in a two dimensional hyperbolic sphere. In the latter case,
we apply the zero-order invariants of the Lie point symmetries in order to
determine invariant solutions of the wave function. Finally, in section \ref%
{conclusions} we draw our conclusions.

\section{Preliminaries}

\label{Preliminaries}

For the convenience of the reader, in this section we discuss briefly the
Lie and Noether point symmetries of differential equations and the
collineations of Riemannian manifolds.

\subsection{Lie point symmetries of differential equations}

Geometrically a differential equation (DE) may be considered as a function $%
H=H(x^{i},u^{A},u_{,i}^{A},u_{,ij}^{A})$ in the space $B=B\left(
x^{i},u^{A},u_{,i}^{A},u_{,ij}^{A}\right) $, where $x^{i}$ are the
independent variables and $u^{A}$ are the dependent variables. The
infinitesimal point transformation
\begin{align}
\bar{x}^{i}& =x^{i}+\varepsilon \xi ^{i}(x^{k},u^{B})~,  \label{pr.01} \\
\bar{u}^{A}& =\bar{u}^{A}+\varepsilon \eta ^{A}(x^{k},u^{B})~,  \label{pr.02}
\end{align}%
has the infinitesimal symmetry generator
\begin{equation}
\mathbf{X}=\xi ^{i}(x^{k},u^{B})\partial _{x^{i}}+\eta
^{A}(x^{k},u^{B})\partial _{u^{A}}~.  \label{pr.03}
\end{equation}

The generator $\mathbf{X}$ of the infinitesimal transformation (\ref{pr.01}%
)-(\ref{pr.02}) is called a Lie point symmetry of the DE $H=0$ if there
exists a function $\kappa $ such that the following condition holds \cite%
{Stephani,Bluman}
\begin{equation}
\mathbf{X}^{[2]}(H)=\kappa H~  \label{pr.04}
\end{equation}%
where
\begin{equation}
\mathbf{X}^{[2]}=\mathbf{X}+\eta _{i}^{A}\partial _{u_{i}^{A}}+\eta
_{ij}^{A}\partial _{u_{ij}^{A}}  \label{pr.05}
\end{equation}%
is the second-prolongation vector of $\mathbf{X,}$ in which
\begin{equation}
\eta _{i}^{A}=\eta _{,i}^{A}+u_{,i}^{B}\eta _{,B}^{A}-\xi
_{,i}^{j}u_{,j}^{A}-u_{,j}^{A}u_{,i}^{B}\xi _{,B}^{j}~,  \label{pr.06}
\end{equation}%
and
\begin{align}
\eta _{ij}^{A}& =\eta _{,ij}^{A}+2\eta _{,B(i}^{A}u_{,j)}^{B}-\xi
_{,ij}^{k}u_{,k}^{A}+\eta _{,BC}^{A}u_{,i}^{B}u_{,j}^{C}-2\xi
_{,(i|B|}^{k}u_{j)}^{B}u_{,k}^{A}  \notag \\
& -\xi _{,BC}^{k}u_{,i}^{B}u_{,j}^{C}u_{,k}^{A}+\eta
_{,B}^{A}u_{,ij}^{B}-2\xi _{,(j}^{k}u_{,i)k}^{A}-\xi _{,B}^{k}\left(
u_{,k}^{A}u_{,ij}^{B}+2u_{(,j}^{B}u_{,i)k}^{A}\right) .  \label{pr.07}
\end{align}

An application of Lie point symmetries of a DE is that they can be used in
order to determine invariant solutions. From the generator (\ref{pr.03}) one
considers the Lagrange system%
\begin{equation}
\frac{dx^{i}}{\xi ^{i}}=\frac{du^{A}}{\eta ^{A}}=\frac{du_{i}^{A}}{\eta _{%
\left[ i\right] }^{A}}=\frac{du_{ij}^{A}}{\eta _{\left[ ij\right] }^{A}}
\end{equation}%
whose solution provides the characteristic functions $W^{\left[ 0\right]
}\left( x^{k},u\right) ,~W^{\left[ 1\right] }\left( x^{k},u,u_{i}\right) $%
~and $W^{\left[ 2\right] }\left( x^{k},u,u_{,i},u_{ij}\right) .$ \ The
characteristic functions can be applied to reduce the order of the DE or the
number of the dependent variables.

Suppose that the DE, $H=H(x^{i},u^{A},u_{,i}^{A},u_{,ij}^{A}),$ arises from
a variational principle, i.e. there exists a Lagrangian function $%
L=L(x^{k},u^{A},u_{,k}^{A})$ such that $H\equiv \mathbf{E}\left( L\right) =0$%
, where $\mathbf{E}$ is the Euler operator. The Lie point symmetry $\mathbf{X%
}$ of the DE\ $H\ \ $is a Noether point symmetry of $H,\ $if the following
additional condition is satisfied%
\begin{equation}
\mathbf{X}^{[1]}L+LD_{i}\xi ^{i}=D_{i}A^{i}\left( x^{k},u^{C}\right) ~
\label{pr.08}
\end{equation}%
where $\mathbf{X}^{[1]}$ is the first prolongation of $\mathbf{X,}$ and $%
D_{i}$ is the covariant derivative wrt the metric $g_{ij}$ in the space of
variables $\{x^{i}\}$, and $A^{i}$ is the Noether current. The
characteristic property of Noether point symmetries is that the quantity
\begin{equation}
I^{i}=\xi ^{k}\left( u_{k}^{A}\frac{\partial L}{\partial u_{i}^{A}}-\delta
_{k}^{i}L\right) -\eta ^{A}\frac{\partial L}{\partial u_{i}^{A}}+A^{i}~
\label{pr.09}
\end{equation}%
is a first integral of Lagrange equations, that is, $D_{i}I^{i}=0,~$\cite%
{Stephani,Bluman}.

\subsection{Collineations of Riemannian spaces}

A collineation in a Riemannian space is a vector field $\mathbf{\xi }$ which
satisfies an equation of the form%
\begin{equation}
\mathcal{L}_{\mathbf{\xi }}\mathbf{A}=\mathbf{B}  \label{pr.10}
\end{equation}%
where $\mathcal{L}_{\xi }$ is the Lie derivative with respect to the vector
field $\mathbf{\xi }$, $\mathbf{A}$ is a geometric object (not necessarily a
tensor) defined in terms of the metric and its derivatives and $\mathbf{B}$
is an arbitrary tensor with the same indices as the geometric object $%
\mathbf{A~}$\cite{Yano}. The collineations of Riemannian spaces have been
classified by Katzin et.al. \cite{Katzin1}. In the following we are
interested in the collineations of the metric tensor\footnote{%
The reason that we consider the symbol $\tilde{g}_{\alpha \beta }\left(
z^{\gamma }\right) $ is because in our case we have two metrics in different
spaces, that is, the metrics $H_{AB}$ and $g_{ij}.$}, i.e. $\mathbf{A}=%
\tilde{g}_{\alpha \beta }\left( z^{\gamma }\right) $, and on the affine
collineations of the connection coefficients $\tilde{\Gamma}_{\beta \gamma
}^{\alpha }\left( z^{\delta }\right) $.

\subsubsection{Conformal symmetries}

The infinitesimal generator $\mathbf{\xi }$ of the point transformation%
\begin{equation}
\bar{z}^{\alpha }=z^{\alpha }+\varepsilon \xi ^{\alpha }\left( z^{\delta
}\right)  \label{pr.11}
\end{equation}%
is called Conformal Killing Vector (CKV) if the Lie derivative of the metric
$\tilde{g}_{\alpha \beta }$ with respect to the vector field $\mathbf{\xi }$
is a multiple of $\tilde{g}_{\alpha \beta }$. That is, if the following
condition holds
\begin{equation}
\mathcal{L}_{\mathbf{\xi }}\tilde{g}_{\alpha \beta }=2\psi \left( z^{\gamma
}\right) \tilde{g}_{\alpha \beta }  \label{pr.12}
\end{equation}%
where $\psi =\frac{1}{n}\xi _{;\gamma }^{\gamma }$.

When $\psi _{;\gamma \delta }=0,$ $\mathbf{\xi }$ is a special\footnote{%
For the conformal factor of a sp.CKV holds $\psi _{;ij}=0,$ that is, $\psi
_{,i}$ is a gradient KV. A Riemannian space admits a sp.CKV if and only if
it admits a gradient KV and a gradient HV \cite{HallspCKV}.} CKV (sp.CKV),
if $\psi =$constant, $\mathbf{\xi }$ is a Homothetic Vector (HV) and when $%
\psi =0$, $\mathbf{\xi }$ is a Killing Vector (KV). \ A metric $\tilde{g}%
_{\alpha \beta }$ admits at most one HV. The CKVs\ of a metric form a Lie
algebra, which is called the conformal algebra, $G_{CV}$. Obviously the KVs
and the homothetic vector are elements of the conformal algebra $G_{CV}$. If
$G_{HV}$ is the algebra of HVs (including the algebra $G_{KV}$ of KVs), then
we have%
\begin{equation}
G_{KV}\subseteq G_{HV}\subseteq G_{CV}  \label{pr.13}
\end{equation}

The maximum dimension of the conformal algebra of an $\ell -~$dimensional
metric $(\ell >2)$ is $G_{\max }=\frac{1}{2}\left( \ell +1\right) \left(
\ell +2\right) ,$ and for $\ell ~=2$ the space admits an infinite
dimensional conformal group. Moreover, if a space $\tilde{g}_{\alpha \beta }$
of dimension $\dim \tilde{g}_{\alpha \beta }>2$ \ admits a conformal algebra
of dimension $G_{\max },$ then it is conformally flat, that is, that there
exists a function $N\left( z^{\gamma }\right) $ such that $\tilde{g}_{\alpha
\beta }=N\left( z^{\gamma }\right) \eta _{\alpha \beta }\,\ $where $\eta
_{\alpha \beta }$ is a flat metric.

CKVs are important in relativistic physics and the effects of the existence
of these vectors can be seen at all levels in General Relativity, that is,
geometry, kinematics and dynamics. We continue with the definition of the
collineations for the connection coefficients of the metric tensor $\bar{g}%
_{\alpha \beta }.$

\subsubsection{Affine collineations}

In a Riemannian space with metric $\tilde{g}_{\alpha \beta }$ and connection
coefficients$~\tilde{\Gamma}_{\beta \gamma }^{\alpha }\left( z^{\delta
}\right) $, the following identity holds
\begin{equation}
\mathcal{L}_{\mathbf{\xi }}\tilde{\Gamma}_{\beta \gamma }^{\alpha }=\tilde{g}%
^{\alpha \delta }\left[ \left( \mathcal{L}_{\mathbf{\xi }}\tilde{g}_{\beta
\delta }\right) _{;\gamma }+\left( \mathcal{L}_{\mathbf{\xi }}\tilde{g}%
_{\delta \gamma }\right) _{;\beta }-\left( \mathcal{L}_{\mathbf{\xi }}\tilde{%
g}_{\beta \gamma }\right) _{;\delta }\right] .  \label{scc.01}
\end{equation}

If $\mathbf{\xi }$ is a HV or KV then from (\ref{scc.01}) follows that $%
\mathcal{L}_{\mathbf{\xi }}\tilde{\Gamma}_{\beta \gamma }^{\alpha }$
vanishes, which implies that the connection coefficients $\tilde{\Gamma}%
_{\beta \gamma }^{\alpha }$ are invariant under the action of transformation
(\ref{pr.11}). In general the infinitesimal generators which leave invariant
the connection coefficients $\tilde{\Gamma}_{\beta \gamma }^{\alpha }$ are
defined by the condition
\begin{equation}
\mathcal{L}_{\mathbf{\xi }}\tilde{\Gamma}_{\beta \gamma }^{\alpha }=0.
\label{scc.02}
\end{equation}%
and are called Affine collineations (AC).

The geometric property of an AC is that it caries a geodesic into a geodesic
and also preserves the affine parameter along each geodesic. \ The ACs of a
Riemannian space form a Lie algebra, which is called the Affine algebra, $%
G_{AC}$ of the space. Obviously the homothetic algebra $G_{HV},~$is a
subalgebra of $G_{AC}$, i.e. $G_{HV}\subseteq G_{AC}.~$We shall say that a
spacetime admits proper ACs when $\dim G_{HV}\prec \dim G_{AC}$. Note that
the proper CKVs do not satisfy condition (\ref{scc.02}) therefore proper
CKVs are not ACs.

In the case of a flat space, condition (\ref{scc.02}) becomes~%
\begin{equation}
\xi _{,\beta \gamma }^{\alpha }=0,  \label{scc.03}
\end{equation}%
whose general solution is~$\xi ^{\alpha }=A_{\beta }^{a}z^{\beta }+B^{\alpha
};~$where $A_{\beta }^{a},B^{\beta }$ are~$\ell \left( \ell +1\right) ~$
constants. Therefore the flat space admits the maximal $\ell \left( \ell
+1\right) $ dimensional Affine algebra. The converse is also true, that is,
if a Riemannian space with metric $\tilde{g}_{\alpha \beta },$ $\dim \tilde{g%
}_{\alpha \beta }=$ $\ell $, admits the Affine algebra $G_{AC}$ with $\dim
G_{AC}=\ell \left( \ell +1\right) $ then the space is flat. We summarize the
above definitions in table \ref{Table11}

%TCIMACRO{\TeXButton{B}{\begin{table}[tcp] \centering}}%
%BeginExpansion
\begin{table}[tcp] \centering%
%EndExpansion
\caption{Collineations of a Riemannian space }%
\begin{tabular}{ccc}
\hline\hline
\textbf{Collineation~}$\mathcal{L}_{\xi }\mathbf{A}=\mathbf{B}$ & \textbf{\ }%
$\mathbf{A}$ & $\mathbf{B}$ \\ \hline
Killing Vector (KV) & $\bar{g}_{ij}$ & $0$ \\
Homothetic vector (HV) & $\bar{g}_{ij}$ & $2\psi g_{ij},~\psi _{,i}=0$ \\
Conformal Killing vector (CKV) & $\bar{g}_{ij}$ & $2\psi g_{ij},~\psi
_{,i}\neq 0$ \\
Affine Collineation (AC) & $\bar{\Gamma}_{jk}^{i}$ & $0$ \\ \hline\hline
\end{tabular}%
\label{Table11}%
%TCIMACRO{\TeXButton{E}{\end{table}}}%
%BeginExpansion
\end{table}%
%EndExpansion

\section{Lie and Noether point symmetries of a class of quasilinear systems
of second-order differential equations}

\label{Lie}

The Lie point symmetry condition (\ref{pr.04}) for the system of equations (%
\ref{pd.02}) has the general form
\begin{equation}
X^{\left[ 2\right] }P^{A}=\kappa _{D}^{A}P^{D}  \label{pd.03}
\end{equation}%
where $\kappa _{D}^{A}$ is a tensor. Replacing (\ref{pr.06}) and (\ref{pr.07}%
), for each term of the left-hand side of condition (\ref{pd.03}) we find
\begin{equation}
\eta ^{D}\frac{\partial P^{A}}{\partial u^{D}}=g^{ij}C_{BC,D}^{A}\eta
^{D}\left( u_{,j}^{B}u_{,i}^{C}\right) +F_{,D}^{A}\eta ^{D}  \label{pd.04}
\end{equation}%
\begin{equation}
\xi ^{k}\frac{\partial P^{A}}{\partial x^{k}}=g_{,k}^{ij}\xi ^{k}\left(
u_{,ij}^{A}\right) +g_{,k}^{ij}\xi ^{k}C_{BC}^{A}\left(
u_{,i}^{B}u_{,j}^{C}\right) -\Gamma _{~,k}^{i}\xi ^{k}\left(
u_{,i}^{A}\right) +F_{,k}^{A}\xi ^{k}  \label{pd.05}
\end{equation}%
\begin{eqnarray}
\eta _{i}^{B}\frac{\partial P^{A}}{\partial u_{,i}^{B}} &=&2g^{ij}C_{BC}^{A}%
\eta _{,i}^{B}\left( u_{,j}^{C}\right) +2g^{ij}C_{BC}^{A}\eta
_{,D}^{B}\left( u_{,j}^{C}u_{,i}^{D}\right) +  \notag \\
&&-2g^{ij}C_{BC}^{A}\xi _{,i}^{k}\left( u_{,k}^{B}u_{,j}^{C}\right)
-2g^{ij}C_{BC}^{A}\xi _{,D}^{k}\left( u_{,i}^{B}u_{,k}^{D}u_{,j}^{C}\right) +
\notag \\
&&-\Gamma ^{i}\eta _{,i}^{A}-\Gamma ^{i}\eta _{,B}^{A}\left(
u_{,i}^{B}\right) +\Gamma ^{i}\xi _{,i}^{j}\left( u_{,j}^{A}\right) +\Gamma
^{i}\xi _{,B}^{j}\left( u_{,i}^{A}u_{,j}^{B}\right) ~  \label{pd.06}
\end{eqnarray}%
and%
\begin{eqnarray}
\eta _{ij}^{B}\frac{\partial P^{A}}{\partial u_{,ij}^{B}} &=&g^{ij}\eta
_{,ij}^{A}+2g^{ij}\eta _{,B(i}^{A}\left( u_{,j)}^{B}\right) -g^{ij}\xi
_{,ij}^{k}\left( u_{,k}^{A}\right) +  \notag \\
&&+g^{ij}\eta _{,BC}^{A}\left( u_{,i}^{B}u_{,j}^{C}\right) -2g^{ij}\xi
_{,(i,|B|}^{k}\left( u_{,j)}^{B}u_{,k}^{A}\right) +  \notag \\
&&-g^{ij}\xi _{,BC}^{k}\left( u_{,i}^{B}u_{,j}^{C}u_{,k}^{A}\right)
+g^{ij}\eta _{,B}^{A}\left( u_{,ij}^{B}\right) +  \notag \\
&&-2g^{ij}\xi _{,(j}^{k}\left( u_{,i)k}^{A}\right) -g^{ij}\xi
_{,B}^{k}\left( u_{,k}^{A}u_{,ij}^{B}+2u_{(,j}^{B}u_{,i)k}^{A}\right) .
\label{pd.07}
\end{eqnarray}%
where indices enclosed in parentheses mean symmetrization, for instance, $%
K_{\left( ij\right) }=\frac{1}{2}\left( K_{ij}+K_{ji}\right) $.

We consider the right-hand side of (\ref{pr.03}) and we introduce new
quantities $\lambda _{D}^{A},\mu _{D}^{k}$ by means of the following relation%
\begin{eqnarray}
\kappa _{D}^{A}P^{D} &=&g^{ij}\lambda _{D}^{A}\left( u_{,ij}^{D}\right)
+g^{ij}\lambda _{D}^{A}C_{BC}^{D}\left( u_{,j}^{B}u_{,i}^{C}\right) -\Gamma
^{i}\lambda _{D}^{A}\left( u_{,i}^{D}\right) +\lambda _{D}^{A}F^{D}+ \\
&&+g^{ij}\mu _{D}^{k}\left( u_{,k}^{A}u_{,ij}^{D}\right) +g^{ij}\mu
_{D}^{k}C_{BC}^{D}\left( u_{,k}^{A}u_{,j}^{B}u_{,i}^{C}\right) -\Gamma
^{i}\mu _{D}^{k}\left( u_{,k}^{A}u_{,i}^{D}\right) +\mu _{D}^{k}F^{D}\left(
u_{,k}^{A}\right) .  \notag
\end{eqnarray}

In order condition (\ref{pd.03}) to hold identically the coefficients of the
terms of the various derivatives of $u^{A}$ in the total expression must be
equal.

For the terms $u_{,k}^{A}u_{(ij)}^{B},~u_{,(i\left\vert k\right\vert
}^{A}u_{j)}^{B},$ we have the following equations

\begin{eqnarray}
u_{,j}^{B}u_{,ik}^{A} &:&\xi _{,B}^{i}=0~,  \label{pd.08} \\
u_{,k}^{A}u_{,ij}^{B} &:&\xi _{,B}^{i}+\mu _{B}^{i}=0~,  \label{pd.08aa}
\end{eqnarray}%
which imply that $\xi ^{i}=\xi ^{i}\left( x^{k}\right) $ and $\mu
_{B}^{k}=0; $ recall that $n\geq 2.$ For $n=1$ we have the case of ODEs for
which the Lie point symmetry conditions are different (see \cite{GNgrg}).

Substituting the solution of the system (\ref{pd.08}), (\ref{pd.08aa}), in (%
\ref{pd.03}) from the coefficients of the remaining terms we obtain the
following symmetry conditions.

Coefficients of $\left( u_{,i}^{A}\right) ^{0}:$%
\begin{equation}
g^{ij}\eta _{,ij}^{A}+F_{,D}^{A}\eta ^{D}+F_{,i}^{A}\xi ^{i}-\Gamma ^{i}\eta
_{,i}^{A}-\lambda _{B}^{A}F^{B}=0.  \label{pd.09}
\end{equation}

Coefficients of $\left( u_{,i}^{A}\right) ^{1}:$%
\begin{eqnarray}
0 &=&-\Gamma _{~,k}^{i}\xi ^{k}\delta _{B}^{A}+\Gamma ^{k}\xi
_{,k}^{i}\delta _{B}^{A}-g^{jk}\xi _{,jk}^{i}\delta _{B}^{A}-\Gamma ^{i}\eta
_{,B}^{A}+  \notag \\
&&+2g^{ik}C_{CB}^{A}\eta _{,k}^{C}+2g^{ik}\eta _{,Bk}^{A}+\lambda
_{B}^{A}\Gamma ^{i}.  \label{pd.10}
\end{eqnarray}

Coefficients of $\left( u_{,i}^{A}\right) ^{2}:$%
\begin{eqnarray}
0 &=&g_{,k}^{ij}\xi ^{k}C_{BC}^{A}+g^{ij}C_{BC,D}^{A}\eta
^{D}+2g^{ij}C_{DB}^{A}\eta _{,C}^{B}+  \notag \\
&&-2g^{ij}C_{BC}^{A}\xi _{,i}^{k}+g^{ij}\eta _{,BC}^{A}-\lambda
_{D}^{A}g^{ij}C_{BC}^{D}.  \label{pd.11}
\end{eqnarray}

Coefficients of$~\left( u_{,ij}^{A}\right) :$%
\begin{equation}
\left( g_{,k}^{ij}\xi ^{k}-2g^{i(k}\xi _{,k}^{j)}\right) \delta
_{B}^{A}+g^{ij}\left( \eta _{,B}^{A}-\lambda _{B}^{A}\right) =0.
\label{pd.12}
\end{equation}

The solution of the system of equations (\ref{pd.09})-(\ref{pd.12}) gives
the generator of the Lie point symmetry vector (\ref{pd.03}). The key point
is to express these conditions in terms of the collineations of the metrics $%
g_{ij,}H_{AB}$ and relate the generator $\mathbf{X}\ $of the Lie point
symmetry to these collineations. Then, in a way, we have geometrized the
problem and we may use the well known results of Differential Geometry in
order to study the Lie point symmetries of the $n\times m$ systems of
differential equations (\ref{pd.02}).

In terms of the Lie derivative equation (\ref{pd.12}) is written as follows,%
\begin{equation}
\left( L_{\xi }g^{ij}\right) \delta _{B}^{A}=-g^{ij}\left( \eta
_{,B}^{A}-\lambda _{B}^{A}\right) .  \label{pd.13}
\end{equation}%
Because $\xi ^{i}=\xi ^{i}\left( x^{k}\right) $ the left-hand side of (\ref%
{pd.13}) is independent of $u^{A}$, hence
\begin{equation}
\lambda _{B}^{A}=\eta _{,B}^{A}-2\psi \left( x^{k}\right) \delta _{B}^{A}.
\label{pd.14}
\end{equation}%
Substituting this back in (\ref{pd.13}) we have
\begin{equation}
\mathcal{L}_{\xi }g_{ij}=2\psi \left( x^{k}\right) g_{ij}  \label{pd.15}
\end{equation}%
which means that $\xi ^{i}\left( x^{k}\right) $ is a CKV of $g_{ij}$ with
conformal factor $\psi \left( x^{k}\right) $. This implies that $\xi
_{;i}^{i}=n\psi \left( x^{k}\right) $, where \textquotedblleft $;$%
\textquotedblright\ indicates covariant derivative with respect to the
metric $g_{ij}$.

Replacing $\lambda _{B}^{A}$ from (\ref{pd.14}) in the symmetry condition (%
\ref{pd.11}) we find
\begin{eqnarray}
0 &=&C_{BC}^{A}\left[ g_{,k}^{ij}\xi ^{k}-2g^{k(j}\xi _{,k}^{i)}+2\psi g^{ij}%
\right] +  \notag \\
&&+g^{ij}\left[ \eta _{,BC}^{A}+C_{BC,D}^{A}\eta ^{D}+2C_{D(B}^{A}\eta
_{,C)}^{D}-\eta _{,D}^{A}C_{BC}^{D}\right]  \notag \\
&=&C_{BC}^{A}\left( \mathcal{L}_{\xi }g^{ij}+2\psi g^{ij}\right) +g^{ij}%
\left[ \eta _{,BC}^{A}+C_{BC,D}^{A}\eta ^{D}+2C_{D(B}^{A}\eta
_{,C)}^{D}-\eta _{,D}^{A}C_{BC}^{D}\right] .  \label{pd.15a}
\end{eqnarray}%
But $C_{BC}^{A}$ are the connection coefficients of the metric $H_{AB}$,
hence
\begin{equation}
\mathcal{L}_{\eta }C_{BC}^{A}=\eta _{,BC}^{A}+C_{BC,D}^{A}\eta
^{D}+2C_{D(B}^{A}\eta _{,C)}^{D}-\eta _{,D}^{A}C_{BC}^{D}.  \label{pd.17}
\end{equation}%
Replacing this back in (\ref{pd.15a}) we find:%
\begin{equation*}
\left( \mathcal{L}_{\xi }g^{ij}+2\psi g^{ij}\right) C_{BC}^{A}+g^{ij}%
\mathcal{L}_{\eta }C_{BC}^{A}=0.
\end{equation*}%
From (\ref{pd.15}) the first term vanishes\footnote{%
Recall that, $\mathcal{L}_{\xi }g^{ij}=-2\psi g^{ij}$.}, therefore the
symmetry condition (\ref{pd.17}) becomes%
\begin{equation}
\mathcal{L}_{\eta }C_{BC}^{A}=0,  \label{pd.18}
\end{equation}%
which means that $\eta ^{A}$ is an AC of $H_{AB}$.

Hence, condition (\ref{pd.10}) becomes%
\begin{equation}
\left( -\Gamma _{~,k}^{i}\xi ^{k}+\Gamma ^{k}\xi _{,k}^{i}-g^{kj}\xi
_{,kj}^{i}-2\psi \Gamma ^{i}\right) \delta _{B}^{A}+2g^{ij}C_{BC}^{A}\eta
_{,j}^{C}+2g^{ij}\eta _{,Bj}^{A}=0.  \label{pd.19}
\end{equation}%
Furthermore, because $\xi ^{i}$ is a CKV of $g_{ij}$ it holds that
\begin{equation}
g^{jk}\mathcal{L}_{\xi }\Gamma _{jk}^{i}=g^{kj}\xi _{,kj}^{i}-\Gamma ^{k}\xi
_{,k}^{i}+\Gamma _{~,k}^{i}\xi ^{k}+2\psi \Gamma ^{i}=\left( 2-n\right) \psi
^{,i}.  \label{pd.20}
\end{equation}

Hence, equation (\ref{pd.19}) becomes%
\begin{equation}
\left( \eta _{|B}^{A}\right) _{,i}=\frac{\left( 2-n\right) }{2}\psi
_{,i}\delta _{B}^{A}  \label{pd.21}
\end{equation}%
where \textquotedblleft $|$\textquotedblright , means covariant derivative
with respect to the metric $H_{AB}$. Furthermore the last equation is
written,%
\begin{equation}
\left( \eta _{A|B}\right) _{,i}=\frac{\left( 2-n\right) }{2}\left( \psi
H_{AB}\right) _{,i}~,
\end{equation}%
from where we have that
\begin{equation}
\eta _{A|B}=\frac{\left( 2-n\right) }{2}\psi H_{AB}+\Lambda _{AB}\left(
u^{C}\right) ,  \label{pd.22}
\end{equation}%
in which, $\Lambda _{AB}=\Lambda _{AB}\left( u^{C}\right) $ is a second rank
tensor defined in the space of the dependent variables. Because $\eta ^{A}$
is an AC of $H_{AB}$ it is true that $\eta _{\left( A|BC\right) }=0$; this
implies that $\Lambda _{\left( AB|C\right) }=0$ which means that $\Lambda
_{AB}$ is a Killing tensor of order two of the metric $H_{AB}$. We conclude
that the general form of $\eta ^{A}\left( x^{k},u^{C}\right) $ is
\begin{equation}
\eta ^{A}=\frac{\left( 2-n\right) }{2}\psi \left( x^{k}\right) Y^{A}\left(
u^{C}\right) +Z^{A}\left( x^{k},u^{B}\right) .  \label{pd.22.1}
\end{equation}

Replacing this in condition (\ref{pd.22}) we find the constraints

\begin{equation}
Y_{A|B}=H_{AB}~~,~~Z_{A|B}=\Lambda _{AB}  \label{pd.23}
\end{equation}%
which mean that vector $Y^{A}$ is a proper gradient HV of $H_{AB}$ and $%
Z^{A} $ is an AC of $H_{AB}$. \ Moreover, when $n>2$ and the space $H_{AB}$
does not admit proper gradient HV then from (\ref{pd.22}) we have that $\psi
\left( x^{k}\right) =0$, i.e. $\xi ^{i}\left( x^{k}\right) $ is a KV of $%
g_{ij}$.

Finally, condition (\ref{pd.09}) gives the further constraint%
\begin{equation}
\left( \mathcal{L}_{\xi }F^{A}+2\psi F^{A}\right) +\frac{2-n}{2}\left( \psi
\mathcal{L}_{Y}F^{A}+g^{ij}\psi _{;ij}Y^{A}\right) +\left( \mathcal{L}%
_{Z}F^{A}+g^{ij}Z_{;ij}^{A}\right) =0.  \label{pd.24}
\end{equation}

The solution of this system which follows from the symmetry condition (\ref%
{pd.03}) leads to the following theorem which is our main result.

\begin{theorem}
\label{Lsym}The Lie point symmetries of the quasilinear systems of
second-order differential equations (\ref{pd.02}) are generated by the CKVs $%
\xi ^{i}\left( x^{k}\right) $ of the metric $g_{ij}$ and the ACs $%
Z^{A}\left( x^{k},u^{C}\right) $ of $\ $the metric $H_{AB}$ such that%
\footnote{%
Where $Z_{A|B}$ means covariant derivative with respect to the metric $%
H_{AB} $ and $\left( Z_{A|B}\right) _{,i}=\frac{\partial }{\partial x^{i}}%
\left( Z_{A|B}\right) $.} $Z_{A|B}=\Lambda _{AB}$ and $\left( Z_{A|B}\right)
_{,i}=0 $ where $\Lambda _{AB}$ is a Killing tensor of order two for the
metric $H_{AB}$ as follows:

(a) If $n>2$ and the metric $H_{AB}$ admits a proper gradient HV $%
~Y^{A}\left( u^{C}\right) ,$ with conformal factor $\bar{\psi}_{Y}=1$, the
generic Lie point symmetry is%
\begin{equation}
X_{L\left( a\right) }=\xi ^{i}\left( x^{k}\right) \partial _{i}+\left[
\left( \frac{2-n}{2}\psi \left( x^{k}\right) \right) Y^{A}\left(
u^{C}\right) +Z^{A}\left( x^{k},u^{C}\right) \right] \partial _{A}
\label{pd.25}
\end{equation}%
and condition (\ref{pd.24}) holds.

(b) If $n>2,$ and the metric $H_{AB}$ does not admit a proper gradient HV,
the generic Lie point symmetry is%
\begin{equation}
X_{L\left( b\right) }=\xi ^{i}\left( x^{k}\right) \partial _{i}+Z^{A}\left(
x^{k},u^{C}\right) \partial _{A}  \label{pd.27}
\end{equation}
and the following condition holds%
\begin{equation}
\mathcal{L}_{\xi }F^{A}+\mathcal{L}_{Z}F^{A}+\Delta _{g}Z^{A}=0.
\label{pd.28}
\end{equation}

(c) If $n=2$, the generic Lie point symmetry is
\begin{equation}
X_{L\left( c\right) }=\xi ^{i}\left( x^{k}\right) \partial _{i}+Z^{A}\left(
x^{k},u^{C}\right) \partial _{A}  \label{pd.28a}
\end{equation}%
and the following condition holds%
\begin{equation}
\left( \mathcal{L}_{\xi }F^{A}+2\psi F^{A}\right) +\left( \mathcal{L}%
_{Z}F^{A}+\Delta _{g}Z^{A}\right) =0.  \label{pd.28b}
\end{equation}
\end{theorem}

We note that Theorem \ref{Lsym} holds for all the systems of the form (\ref%
{pd.02}) i.e. they do not necessarily admit a Lagrangian.

\subsection{Noether symmetries}

In this section we study the Noether point symmetries of Lagrangian (\ref%
{pd.01}). For each term of the Noether condition (\ref{pr.08}) for the
Lagrangian (\ref{pd.01}) we have

\begin{equation}
D_{i}A^{i}=A_{,i}^{i}+A_{,A}^{i}\left( u_{,i}^{A}\right) ,  \label{pd.29}
\end{equation}%
\begin{eqnarray}
LD_{i}\xi ^{i} &=&\left( \frac{1}{2}\sqrt{g}g^{ij}H_{AB}\xi _{,k}^{k}\right)
\left( u_{,i}^{A}u_{,j}^{B}\right) -\sqrt{g}V\xi _{,k}^{k}+  \notag \\
&&+\left( \frac{1}{2}\sqrt{g}g^{ij}H_{AB}\xi _{,C}^{k}\right) \left(
u_{,i}^{A}u_{,j}^{B}u_{,k}^{C}\right) -\sqrt{g}\xi _{,A}^{k}V~\left(
u_{,k}^{A}\right) .  \label{pd.30}
\end{eqnarray}%
Moreover, for the terms of $X^{\left[ 1\right] }L$ we find,%
\begin{equation}
\eta ^{C}\frac{\partial L}{\partial u^{C}}=\left( \frac{1}{2}\sqrt{g}%
g^{ij}\eta ^{C}H_{AB,C}\right) \left( u_{,i}^{A}u_{,j}^{B}\right) -\sqrt{g}%
V_{,C}\eta ^{C},  \label{pd.31}
\end{equation}%
\begin{equation}
\xi ^{k}\frac{\partial L}{\partial x^{k}}=\left( \frac{1}{2}\xi ^{k}\left(
\sqrt{g}g^{ij}\right) _{,k}H_{AB}\right) \left( u_{,i}^{A}u_{,j}^{B}\right)
-\xi ^{k}\left( \sqrt{g}V\right) _{,k},  \label{pd.32}
\end{equation}%
and%
\begin{eqnarray}
\eta _{k}^{C}\frac{\partial L}{\partial u_{,k}^{C}} &=&\left( \sqrt{g}%
g^{ij}H_{AB}\eta _{,i}^{A}\right) \left( u_{,j}^{B}\right) -\left( \sqrt{g}%
g^{ij}H_{AB}\xi _{,C}^{k}\right) \left( u_{,j}^{B}u_{,i}^{A}u_{,k}^{C}\right)
\notag \\
&&+\left( \sqrt{g}g^{ij}H_{AB}\eta _{,C}^{A}-\sqrt{g}g^{ij}H_{AB}\xi
_{,i}^{j}\right) \left( u_{,j}^{B}u_{,i}^{C}\right) .  \label{pd.33}
\end{eqnarray}

From the coefficients of the monomial $\left( u_{,i}^{A}\right) ^{3}$ we
have $\xi _{,C}^{j}=0$, i.e. $\xi ^{i}=\xi ^{i}\left( x^{k}\right) $. This
should be expected because the Noether point symmetries are Lie point
symmetries for which (as we have shown already) $\xi ^{i}=\xi ^{i}\left(
x^{k}\right) .$

Replacing (\ref{pd.29})-(\ref{pd.33}) in the Noether condition (\ref{pr.08})
and using the Lie derivative we find the following Noether symmetry
conditions.

Coefficients of $\left( u_{,i}^{A}\right) ^{0}:$%
\begin{equation}
\sqrt{g}\left( \mathcal{L}_{\eta }V+\mathcal{L}_{\xi }V+\xi
_{;k}^{k}V\right) +A_{,k}^{k}=0  \label{pd.35}
\end{equation}

Coefficients of $\left( u_{,i}^{A}\right) ^{1}:$%
\begin{equation}
\sqrt{g}g^{ij}H_{AB}\eta _{,j}^{A}-A_{,B}^{i}=0  \label{pd.34}
\end{equation}

Coefficients of $\left( u_{,i}^{A}\right) ^{2}:$%
\begin{equation}
H_{AB}\left( \mathcal{L}_{\xi }g^{ij}+\xi _{,k}^{k}g^{ij}\right)
+g^{ij}\left( \mathcal{L}_{\eta }H_{AB}\right) =0.  \label{pd.36}
\end{equation}

From Theorem \ref{Lsym}, we know that $\xi ^{i}$ is a CKV of $g_{ij}$; hence
$\xi _{;k}^{k}=n\psi \left( x^{k}\right) $. Substituting in (\ref{pd.36}) we
find%
\begin{equation}
H_{AB}\left( \mathcal{L}_{\xi }g_{ij}\right) =g_{ij}\left( \mathcal{L}_{\eta
}H_{AB}+n\psi H_{AB}\right)  \label{pd.37}
\end{equation}%
or equivalently%
\begin{equation}
\mathcal{L}_{\eta }H_{AB}=\left( 2-n\right) \psi H_{AB}  \label{pd.38}
\end{equation}%
which implies that
\begin{equation}
\eta ^{A}\left( x^{k},u^{C}\right) =\frac{2-n}{2}\psi \left( x^{k}\right)
Y^{A}\left( u^{K}\right) +K^{A}\left( x^{k},u^{C}\right)  \label{pd.39}
\end{equation}%
where $Y^{A}\left( u^{k}\right) $ is a proper gradient HV of $H_{AB}$ and $%
K^{A}\left( x^{k},u^{C}\right) $ is a KV of $H_{AB}$.

Substituting back in (\ref{pd.34}) we have
\begin{equation}
A_{i,A}=\frac{2-n}{2}\sqrt{g}\psi _{,i}Y_{A}+K_{A,i}  \label{pd.40}
\end{equation}%
which gives that $Y_{A}$ is a gradient HV of $H_{AB}$, i.e. $Y_{A}=Y_{,A}$
and $K^{A}=K\left( x^{k},u^{C}\right) $. Moreover if $K^{A}$ is a non
gradient KV of $H_{AB}$ then from (\ref{pd.40}) we have that $%
K^{A}=K^{A}\left( u^{C}\right) $.

Therefore, we may write $K^{A}=K_{G}^{A}(x^{k},u^{C})+K_{NG}^{A}\left(
u^{C}\right) $ where $K_{G}^{A},~K_{NG}^{A}$ are the gradient and non
gradient KVs respectively. Then from (\ref{pd.40}) we have for following
expression for the Noether vector%
\begin{equation}
A_{i}=\frac{2-n}{2}\sqrt{g}\psi _{,i}Y+\sqrt{g}K_{G,i}+\sqrt{g}\Phi
_{i}\left( x^{k}\right) .  \label{pd.41}
\end{equation}%
Finally from (\ref{pd.36}) we have the constraint%
\begin{equation}
\mathcal{L}_{\xi }V+n\psi V+\frac{2-n}{2}\left( \psi V_{,A}Y^{A}+\Delta
_{g}\psi Y\right) +\left( \Delta _{g}\left( K_{G}\right) +V_{,A}K^{A}\right)
+\Phi _{,k}^{k}=0.  \label{pd.42}
\end{equation}

We collect the results in the following theorem.

\begin{theorem}
\label{Nsym}The Noether point symmetries of Lagrangian (\ref{pd.01}) are
generated by the CKVs $\xi ^{i}$ \ of the metric $g_{ij}$ with conformal
factor $\psi \left( x^{k}\right) $ and the KVs of $\ $the metric $H_{AB}$ as
follows:

(a) If $n>2$ and the metric $H_{AB}$ admits a proper gradient HV with
homothetic factor $\bar{\psi}_{Y}=1$, the generic Noether point symmetry is%
\begin{equation}
X_{N\left( a\right) }=\xi ^{i}\left( x^{k}\right) \partial _{i}+\left[
\left( \frac{2-n}{2}\psi \left( x^{k}\right) \right) Y^{A}\left(
u^{C}\right) +K_{G}^{A}(x^{k},u^{C})+K_{NG}^{A}\left( u^{C}\right) \right]
\partial _{A}  \label{pd.43}
\end{equation}%
where , $K_{G}^{A}$ is a gradient KV/HV of $H_{AB}$ , $K_{NG}^{A}$ is a non
gradient KV/HV of $H_{AB}$ and condition (\ref{pd.42}) holds.

The corresponding gauge vector field is
\begin{equation}
A_{\left( a\right) }^{i}\left( x^{k},u^{C}\right) =\sqrt{g}\left( \frac{2-n}{%
2}\sqrt{g}g^{ij}\psi _{,j}Y+g^{ij}K_{G,j}+\Phi ^{i}\left( x^{k}\right)
\right)  \label{pd.45}
\end{equation}%
and the generic Noether conservation current is
\begin{eqnarray}
I_{\left( a\right) }^{i} &=&\xi ^{k}\mathcal{H}_{k}^{i}-\left( \left( \frac{%
2-n}{2}\psi \right) Y^{A}+K_{G}^{A}+K_{NG}^{A}\right) g^{ij}H_{AB}u_{,j}^{A}+
\notag \\
&&+\sqrt{g}\left( \frac{2-n}{2}\sqrt{g}g^{ij}\psi _{,j}Y+g^{ij}K_{G,j}+\Phi
^{i}\right) .  \label{pd.45a}
\end{eqnarray}

(b) If $n>2,$ and the metric $H_{AB}$ does not admit a proper gradient HV,
the generic Noether point symmetry is
\begin{equation}
X_{N\left( b\right) }=\xi ^{i}\left( x^{k}\right) \partial _{i}+\left(
K_{G}^{A}(x^{k},u^{C})+K_{NG}^{A}\left( u^{C}\right) \right) \partial _{A}
\label{pd.46}
\end{equation}%
where $K_{G}^{A}$ is a gradient KV of $H_{AB}$ , $K_{NG}^{A}$ is a non
gradient KV of $H_{AB}$ and the following condition holds%
\begin{equation}
\mathcal{L}_{\xi }V+\left( \Delta _{g}\left( K_{G}\right)
+V_{,A}K^{A}\right) +\Phi _{,k}^{k}=0.  \label{pd.47}
\end{equation}%
The corresponding gauge vector field is
\begin{equation}
A_{\left( b\right) }^{i}\left( x^{k},u^{C}\right) =\sqrt{g}\left(
g^{ij}K_{G,j}+\Phi ^{i}\left( x^{k}\right) \right)  \label{pd.48}
\end{equation}%
and the generic Noether conservation current is
\begin{equation}
I_{\left( b\right) }^{i}=\xi ^{k}\mathcal{H}_{k}^{i}-\left(
K_{G}^{A}+K_{NG}^{A}\right) g^{ij}H_{AB}u_{,j}^{B}+\sqrt{g}\left(
g^{ij}K_{G,j}+\Phi ^{i}\right)  \label{pd.48b}
\end{equation}

(c) If $n=2$, the generic Noether point symmetry is
\begin{equation}
X_{N\left( c\right) }=\xi ^{i}\left( x^{k}\right) \partial _{i}+\left(
K_{G}^{A}(x^{k},u^{C})+K_{NG}^{A}\left( u^{C}\right) \right) \partial _{A}
\label{pd.49}
\end{equation}%
where $K_{G}^{A}$ is a gradient KV/HV of $H_{AB}$ , $K_{NG}^{A}$ is a non
gradient KV/HV of $H_{AB}$ and the following condition holds%
\begin{equation}
\mathcal{L}_{\xi }V+2\psi V+\left( \Delta _{g}\left( K_{G}\right)
+V_{,A}K^{A}\right) +\Phi _{,k}^{k}=0.  \label{pd.50}
\end{equation}%
The corresponding gauge vector field is%
\begin{equation}
A_{\left( c\right) }^{i}\left( x^{k},u^{C}\right) =\sqrt{g}\left(
g^{ij}K_{G,j}+\Phi ^{i}\left( x^{k}\right) \right)  \label{pd.51}
\end{equation}%
and the generic Noether conservation current is
\begin{equation}
I_{\left( c\right) }^{i}=\xi ^{k}\mathcal{H}_{k}^{i}-\left(
K_{G}^{A}+K_{NG}^{A}\right) g^{ij}H_{AB}u_{,j}^{B}+\sqrt{g}\left(
g^{ij}K_{G,j}+\Phi ^{i}\right) .  \label{pd.51b}
\end{equation}%
\newline
In all cases the function $\mathcal{H}=\mathcal{H}\left(
x^{k},u^{C},u_{,k}^{C}\right) $, is the Hamiltonian of Lagrangian (\ref%
{pd.01}); that is,%
\begin{equation}
\mathcal{H}_{k}^{i}=\frac{1}{2}\sqrt{g}H_{AB}\left(
2g^{ij}u_{k}^{(A}u_{,j}^{B)}-\delta
_{k}^{i}g^{rs}u_{,r}^{A}u_{,s}^{B}\right) +\delta _{k}^{i}\sqrt{g}V
\end{equation}
\end{theorem}

The vector fields $X_{N\left( a\right) },~X_{N\left( b\right) }~$and $%
X_{N\left( c\right) }$ of Theorem \ref{Nsym} give the generic Noether point
symmetry of Lagrangian (\ref{pd.01}) in a Riemannian space $g_{ij}$.

In the following sections we proceed with the applications of Theorems \ref%
{Lsym} and \ref{Nsym} in two cases of special interest. Specifically, we
study the point symmetries of a system of quasilinear Laplace equations, and
the point symmetries of the modified Klein-Gordon equation for a particle in
Generalized Uncertainty Principle.

\section{System of quasilinear Laplace equations}

\label{Laplace}

We assume that the potential $V\left( x^{k},u^{C}\right) $ of (\ref{pd.01})
is zero. Then the Euler-Lagrange equations (\ref{pd.02}) become
\begin{equation}
g^{ij}u_{,ij}^{A}+g^{ij}C_{BC}^{A}u_{,j}^{B}u_{,i}^{C}-\Gamma
^{i}u_{,i}^{A}=0  \label{pd.53}
\end{equation}%
and correspond to a system of quasilinear Laplace of dimension $m$. When $%
n=1 $ the system (\ref{pd.53}) describes the geodesic equations of a
particle with affine parameterization in the space $H_{AB}$, or the wave
equation in the space $g_{ij}$ when $m=1$. Recall that in this work we
consider $n\geq 2$. For $n=1$ we have the case of autoparallel equations see
\cite{GRG1}.

For this particular case, from Theorems \ref{Lsym} and \ref{Nsym} we have
the following corollary.

\begin{corollary}
\label{CorA}The generic form of the generators of Lie and Noether point
symmetries of the system of second-order PDEs (\ref{pd.53}) are those of
Theorems \ref{Lsym} and \ref{Nsym}, where the corresponding constraint
conditions are as follows:

(a) If $n\geq 2$ and the metric $H_{AB}$ admits a proper gradient HV, the
Lie point symmetry constraint condition is $\frac{2-n}{2}\Delta _{g}\psi
Y^{A}+\Delta _{g}Z^{A}=0~\ $and the Noether point symmetry condition becomes
$\frac{2-n}{2}\Delta _{g}\psi Y+\Delta _{g}\left( K_{G}\right) +\Phi
_{,k}^{k}=0$.

(b,c) If $n>2,$ and the metric $H_{AB}$ does not admit a proper gradient HV,
or if $\dim g_{ij}=2$, the Lie constraint condition is $\Delta _{g}Z^{A}=0$
and the Noether symmetry condition becomes $\Delta _{g}K_{G}=0$.
\end{corollary}

We observe that in this particular application the main role is played by
the metric $H_{AB}\left( u^{C}\right) .$ Therefore, we study two important
cases; that is, (a) the space $H_{AB}\left( u^{C}\right) $ is flat, \ (b)
and $H_{AB}\left( u^{C}\right) $ is a space of constant curvature.

\subsection{Case a: $H_{AB}$ is flat}

We consider a Euclidean space of dimension $m$ in which we employ Cartesian
coordinates so that $H_{AB}=\delta _{AB}$. In this case the system (\ref%
{pd.53}) takes the simplest form:%
\begin{equation}
g^{ij}u_{,ij}^{A}-\Gamma ^{i}u_{,i}^{A}=0.  \label{pd.53a}
\end{equation}

\ As we have already remarked, the $m$ dimensional flat space admits an $%
m\left( m+1\right) $ dimensional Lie algebra of ACs. This algebra consists
of $m$ linearly independent gradient KVs and $m^{2}$ proper ACs. We note
that the gradient HV and the non gradient KVs (rotation group) of the flat
space follow from linear combinations of the proper ACs. In table \ref%
{ACsFlat} we give the KVs, the HV and the proper ACs of the one, two and
three dimensional flat space.
%TCIMACRO{\TeXButton{B}{\begin{table}[tbp] \centering}}%
%BeginExpansion
\begin{table}[tbp] \centering%
%EndExpansion
\caption{Affine collineations for flat space of dimension $m$, with
$m=1,2,3$}%
\begin{tabular}{ccc}
\hline\hline
$\mathbf{m=1}$ & \textbf{Gradient} & \textbf{Non gradient} \\
\textbf{KV} & $\partial _{u^{1}}$ & $\nexists $ \\
\textbf{HV} & $u^{1}\partial _{u^{1}}$ & $\nexists $ \\
\textbf{ACs} & $\nexists $ & $\nexists $ \\
$\mathbf{m=2}$ & \textbf{Gradient} & \textbf{Non gradient} \\ \hline
\textbf{KVs} & $\partial _{u^{1}}~,~\partial _{u^{2}}$ & $u^{2}\partial
_{u^{1}}-u^{1}\partial _{u^{2}}$ \\
\textbf{HV} & $u^{1}\partial _{u^{1}}+u^{2}\partial _{u^{2}}$ & $\nexists $
\\
\textbf{ACs} & $u^{1}\partial _{u^{1}}~,~u^{2}\partial _{u^{2}}$ & $%
u^{2}\partial _{u^{1}}~,~u^{1}\partial _{u^{2}}$ \\
$\mathbf{m=3}$ &  &  \\ \hline
\textbf{KVs} & $\partial _{u^{1}}~,~\partial _{u^{2}}~,~\partial _{u^{3}}$ &
$u^{2}\partial _{u^{1}}-u^{1}\partial _{u^{2}}~,~u^{3}\partial
_{u^{1}}-u^{1}\partial _{u^{3}}~$ \\
&  & $u^{3}\partial _{u^{2}}-u^{2}\partial _{u^{3}}$ \\
\textbf{HV} & $u^{1}\partial _{u^{1}}+u^{2}\partial _{u^{2}}+u^{3}\partial
_{u^{3}}$ & $\nexists $ \\
\textbf{ACs} & $u^{1}\partial _{u^{1}}~,~u^{2}\partial
_{u^{2}}~,~u^{3}\partial _{u^{3}}$ & $u^{2}\partial _{u^{1}}~,~u^{3}\partial
_{u^{1}}~,~u^{1}\partial _{u^{2}}$ \\
&  & $u^{3}\partial _{u^{2}}~,~u^{1}\partial _{u^{3}}~,~u^{2}\partial
_{u^{3}}$ \\ \hline\hline
\end{tabular}%
\label{ACsFlat}%
%TCIMACRO{\TeXButton{E}{\end{table}}}%
%BeginExpansion
\end{table}%
%EndExpansion

The general AC for the $m$ dimensional flat space is%
\begin{equation}
Z^{A}=b^{A}\partial _{u^{A}}+c_{B}^{A}u^{B}\partial _{u^{A}}  \label{pd.54}
\end{equation}%
where $b^{A},c_{B}^{A}$ are constants in the space $H_{AB}$ so that $%
b^{A}=b^{A}\left( x^{k}\right) $ and $c_{B}^{A}=c_{B}^{A}\left( x^{k}\right)
.$ Moreover from the condition $\left( Z_{A|B}\right) _{,i}=0~$of Theorem %
\ref{Lsym}, we have that $c_{B,i}^{A}=0$ and $b^{A}=b^{A}\left( x^{k}\right)
$. Then from Corollary \ref{CorA}\ follows that the generic Lie symmetry
vector of the system (\ref{pd.53a}) is
\begin{equation}
X_{L}=\xi ^{i}\left( x^{k}\right) \partial _{i}+\left( \frac{2-n}{2}\psi
\left( x^{k}\right) \right) u^{A}\partial _{A}+b^{A}\left( x^{k}\right)
\partial _{A}+c_{B}^{A}u^{B}\partial _{u^{A}}  \label{pd.55}
\end{equation}%
where the following condition holds%
\begin{equation}
\frac{2-n}{2}\Delta _{g}\psi =0~,~\Delta _{g}b^{A}=0.  \label{pd.56}
\end{equation}

From the second condition it follows that the functions $b^{A}\left(
x^{k}\right) $ are solutions of (\ref{pd.53a}) and the conformal factor $%
\psi \left( x^{k}\right) $ satisfies the Laplacian in the space of the
independent variables with metric $g_{ij}$.

The generic Noether symmetry vector for the Lagrangian (\ref{pd.01}) for the
system (\ref{pd.53a}) is the vector field (\ref{pd.55}) with the constraints
(\ref{pd.56}), where now the constants $c_{AB}=c_{IJ}\delta _{A}^{[I}\delta
_{B}^{J]}~$with $c_{IJ}\in
%TCIMACRO{\U{211d} }%
%BeginExpansion
\mathbb{R}
%EndExpansion
$. That means that the components $c_{B}^{A}u^{B}\partial _{u^{A}}$ of (\ref%
{pd.55}) are HVs of $\delta _{AB}$.

When $\dim H_{AB}=1$, from Table \ref{ACsFlat}, it follows that the space
admits a two dimensional affine algebra; the two vector fields are a
gradient KV and a HV. Therefore, what it has been called as "linear/trivial"
symmetries of Laplace equation (with $m=1$) \cite{AnIJGMMP,Bozhkov},$\ $are
the symmetries which arise from the gradient KVs/HV of the one dimensional
space.

Furthermore, we consider $g_{ij}$ to be the flat space metric with dimension
$n>2,~$i.e. $g_{ij}=\delta _{ij}$. Then the system (\ref{pd.53a}) takes the
simplest form $\delta ^{ij}u_{,ij}^{A}=0$ which corresponds to a system of $%
m-$Laplace equations. It is well known that the flat space $\delta _{ij}$
admits a $\frac{\left( n+2\right) \left( n+1\right) }{2}$ dimensional
conformal algebra where the proper CKVs are $n$, with the property $\psi
_{,ij}=0$, i.e. they are gradient. Furthermore, the case where the two
metrics $g_{ij}$, $H_{AB}$ are flat corresponds to the case in which
equation (\ref{pd.53a}) admits the maximum Lie point symmetries. We conclude
with the following corollary.

\begin{corollary}
\label{CorCounting} Consider the $n\times m$ system of second-order system
of PDEs (\ref{pd.53a}) with $n>2$ and Lagrangian (\ref{pd.01}). Then:\newline
(A) If the system (\ref{pd.53a}) is invariant under the action of the group $%
\tilde{G}~$(Lie point symmetries)$,$ then $\dim \tilde{G}\leq \frac{\left(
n+2\right) \left( n+1\right) }{2}+m\left( m+1\right) $.\newline
(B) If the Lagrangian (\ref{pd.01}) with $V\left( x^{k},u^{C}\right) =0$,
admits Noether point symmetries which form the group $\tilde{G}_{N}$, then $%
\dim \tilde{G}_{N}$ $\leq \frac{\left( n+2\right) \left( n+1\right) }{2}+%
\frac{m\left( m+1\right) }{2}$.\newline
In both cases (A)\ and (B)\ the equality holds when and only when $g_{ij},$ $%
H_{AB}$ are flat spaces of dimension $n,$ and $m$ respectively. In this case
there exists a coordinate system such that the system (\ref{pd.53a}) becomes
$\delta ^{ij}u_{,ij}^{A}=0$.
\end{corollary}

We would like to remark, that the results of this subsection hold and in the
case the flat spaces $H_{AB}$ or $g_{ij}$ have Lorentzian signature. What
changes in this case is the form of the non gradient KVs of Table \ref%
{ACsFlat}. Additionally, for $m=1$ Theorem \ref{CorCounting} gives the
results for the wave equation \cite{ibrawave}. Moreover, in the case of the
geodesic equations, i.e. $n=1$, the maximum Noether algebra is consistent
with that of the geodesic Lagrangian (\ref{pd.01}). However, when $n=1$, the
maximum algebra of Lie point symmetries is different from that of theorem %
\ref{CorCounting}, and it is $\left( m+1\right) \left( m+3\right) $, which
is the projective algebra of the $m+1$ flat spacetime \cite{Aminova,Aminova2}%
. That is, in the case of geodesic equations the Lie point symmetries follow
from the special projective algebra of the space $H_{AB}$.

\subsection{Case b: $H_{AB}$ is the metric of a space of constant curvature}

We assume now that $H_{AB}$ is a space of constant curvature$\ K$ with~$%
K\neq 0.$ We choose coordinates so that $H_{AB}=U\left( u^{C}u_{C}\right)
\delta _{AB},$ and $U\left( u^{C}\right) =\left( 1+\frac{K}{4}\delta
_{AB}u^{A}u^{B}\right) ^{-2};$ in field theory the models which live in that
space are called $\sigma -$models \cite{Chervon1}.

In these coordinates, the connection coefficients of the metric $H_{AB}$
have the following form%
\begin{equation}
C_{BC}^{A}=-\frac{KU}{2}\left( u_{C}\delta _{B}^{A}+u_{B}\delta
_{C}^{A}-u^{A}\delta _{BC}\right)  \label{pd.57}
\end{equation}%
and equation (\ref{pd.53a}) becomes%
\begin{equation}
g^{ij}u_{,ij}^{A}-\frac{KU}{2}g^{ij}u_{,j}^{B}u_{,i}^{C}\left( u_{C}\delta
_{B}^{A}+u_{B}\delta _{C}^{A}-u^{A}\delta _{BC}\right) -\Gamma
^{i}u_{,i}^{A}=0.  \label{pd.58}
\end{equation}

In order to determine the Lie and the Noether point symmetries of the system
(\ref{pd.58}) we have to study the Affine algebra of a space of constant
curvature. It is well known that the Affine algebra of a \ space of constant
non-vanishing curvature is the $SO\left( n+1\right) $ Lie algebra of non
gradient KVs whose dimension is $\dim SO\left( n+1\right) =\frac{n\left(
n+1\right) }{2}$ \cite{Barnes}. Therefore the Lie and the Noether point
symmetries of (\ref{pd.58}) follow from theorem \ref{CorA}(b). We have the
following result.

\begin{corollary}
\textbf{Theorem} \label{CorCountingB} If the $n\times m$ system of
second-order PDEs (\ref{pd.58}) with $n>2$ is invariant under the action of
the group $\hat{G}$, then $\dim \hat{G}\leq \frac{n\left( n+1\right) }{2}+%
\frac{m\left( m-1\right) }{2}$. The equality holds when and only when $%
g_{ij} $ is the metric of a maximally symmetric space. Moreover, all the Lie
point symmetries of (\ref{pd.58}) are also Noether point symmetries of
Lagrangian (\ref{pd.01}) with $V\left( x^{k},u^{C}\right) =0$.
\end{corollary}

\section{The Klein-Gordon equation modified by the Generalized Uncertainty
Principle}

\label{GUPkg} In this section, we study the Lie and the Noether point
symmetries of the modified Klein-Gordon equation of a spin-0 particle
modified by the Generalized Uncertainty Principle (GUP) \cite%
{Maggiore,Kemph1,Vagenas,Moayedi}. The modified Klein-Gordon equation is a
fourth order PDE which, by means of a Lagrange multiplier, is reduced to a
system of two second-order PDEs of the form of the system (\ref{pd.02}). By
applying the results of sections \ref{Lsym} and \ref{Nsym}, we prove a
corollary concerning the modified Klein-Gordon equation in a Riemannian
space $g_{ij}$. We apply the corollary in two cases of physical interest:
(A) the underlying manifold is the flat space $M^{4}$, and (B) the
underlying manifold is the two dimensional hyperbolic sphere; in each case
we determine the Lie and the Noether point symmetries, and we apply the
zero-order invariants in order to determine invariant solutions of the wave
function.

\subsection{Generalized Uncertainty Principle}

The modified structural form of GUP is
\begin{equation}
\Delta X_{i}\Delta P_{j}\geqslant \frac{\hbar }{2}[\delta _{\alpha \beta
}(1+\beta P^{2})+2\beta P_{\alpha }P_{\beta }]  \label{pd.59}
\end{equation}%
where the deformed Heisenberg algebra which is found from (\ref{pd.59}) is
\begin{equation}
\lbrack X_{i},P_{j}]=i\hbar \lbrack \delta _{\alpha \beta }(1+\beta
P^{2})+2\beta P_{\alpha }P_{\beta }].  \label{pd.60}
\end{equation}%
Here $\beta $ is a parameter of deformation defined by\footnote{$M_{Pl}$ is
the Planck mass, $\ell _{Pl}$ $(\approx 10^{-35}~m)$ is the Planck length, $%
M_{Pl}c^{2}$ $(\approx 1.2~10^{19}~GeV)$ is the Planck energy.} $\beta ={%
\beta _{0}}/{M_{Pl}^{2}c^{2}}={\beta _{0}\ell _{Pl}^{2}}/{\hbar ^{2}}$. By
keeping $X_{\alpha }=x_{\beta }$ undeformed, the coordinate representation
of the momentum operator is $P_{\alpha }=p_{\alpha }(1+\beta p^{2})$; $%
\left( x,p\right) $ is the canonical representation satisfying $[x_{\alpha
},p_{\beta }]=i\hbar \delta _{\alpha \beta }$.

In the relativistic four vector form, the commutation relation (\ref{pd.60})
can be written as \cite{Moayedi}
\begin{equation}
\lbrack X_{\mu },P_{\nu }]=-i\hbar \lbrack (1-\beta (\eta ^{\mu \nu }P_{\mu
}P_{\nu }))\eta _{\mu \nu }-2\beta P_{\mu }P_{\nu }]  \label{pd.61}
\end{equation}%
where $\eta _{\mu \nu }=diag(1,-1,-1,-1)$.{\LARGE \ }The corresponding
deformed operators in this case are
\begin{equation}
P_{\mu }=p_{\mu }(1-\beta (\eta ^{\alpha \gamma }p_{\alpha }p_{\gamma
}))~,~~X_{\nu }=x_{\nu },  \label{pd.62}
\end{equation}%
where$~p^{\mu }=i\hbar \frac{\partial }{\partial x_{\mu }},$ and $[x_{\mu
},p_{\nu }]=-i\hbar \eta _{\mu \nu }$.

Consider a spin-0 particle with rest mass $m.~$The Klein-Gordon equation of
this particle is%
\begin{equation}
\left[ \eta ^{\mu \nu }P_{\mu }P_{\nu }-\left( mc\right) ^{2}\right] \Psi =0
\end{equation}%
where $c$ is the speed of light. By substituting $P_{\mu }$ from (\ref{pd.62}%
), we have the modified Klein-Gordon equation%
\begin{equation}
\Delta \Psi -2\beta \hbar ^{2}\Delta \left( \Delta \Psi \right) +V_{0}\Psi
=0,  \label{pd.63}
\end{equation}%
where $V_{0}=\left( \frac{mc}{\hbar }\right) ^{2};$ $\Delta $ is the Laplace
operator where in $M^{4}$, $\Delta \equiv \square ~,$ and the terms $O\left(
\beta ^{2}\right) $ have been eliminated. Equation (\ref{pd.63}) is a fourth
order PDE; however, with the use of a Lagrange multiplier we can write it as
a system of two second-order PDEs.

The action of the modified Klein-Gordon equation (\ref{pd.63}) is%
\begin{equation}
S=\int dx^{4}\sqrt{-g}L\left( \Psi ,\mathcal{D}_{\sigma }\Psi \right)
\label{pd.64}
\end{equation}%
where the Lagrangian $L\left( \Psi ,\mathcal{D}_{\sigma }\Psi \right) $ of
the models is

\begin{equation}
L\left( \Psi ,\mathcal{D}_{\sigma }\Psi \right) =\frac{1}{2}\sqrt{-g}g^{\mu
\nu }\mathcal{D}_{\mu }\Psi \mathcal{D}_{\nu }\Psi -\frac{1}{2}\sqrt{-g}%
V_{0}\Psi ^{2}  \label{pd.65}
\end{equation}%
and the new operator $\mathcal{D}_{\mu }$ is $\mathcal{D}_{\mu }=\nabla
_{\mu }+\beta \hbar ^{2}\nabla _{\mu }\left( \Delta \right) $;$~\nabla _{\mu
}$ is the covariant derivative, i.e.$~\nabla _{\mu }\Psi =\Psi _{;\mu }$. We
introduce the new variable $\Phi =\Delta _{g}\Psi $, and the Lagrange
multiplier $\lambda .$ From the constraint $\frac{\delta S}{\delta \lambda }%
=0$ we have that $\lambda =-2\beta \hbar ^{2}\Phi $, and the action (\ref%
{pd.64}) becomes%
\begin{equation}
S=\int dx^{4}\sqrt{-g}\left( \frac{1}{2}g^{\mu \nu }\Psi _{;\mu }\Psi _{;\nu
}+2\beta \hbar ^{2}g^{\mu \nu }\Psi _{;\mu }\Phi _{;\nu }+\beta \hbar
^{2}\Phi ^{2}-\frac{1}{2}V_{0}\Psi ^{2}\right) .  \label{pd.66}
\end{equation}

Hence, the new Lagrangian is%
\begin{equation}
L\left( \Psi ,\Psi _{;\mu },\Phi ,\Phi _{;\mu }\right) =\sqrt{-g}\left(
\frac{1}{2}g^{\mu \nu }\Psi _{;\mu }\Psi _{;\nu }+2\beta \hbar ^{2}g^{\mu
\nu }\Psi _{;\mu }\Phi _{;\nu }\right) -\sqrt{-g}\left( \frac{1}{2}V_{0}\Psi
^{2}-\beta \hbar ^{2}\Phi ^{2}\right)  \label{pd.67}
\end{equation}%
where $\Phi $ is a new field. We note that the Lagrangian (\ref{pd.67}) is
of the form (\ref{pd.01}) with%
\begin{equation}
H_{AB}=%
\begin{pmatrix}
1 & 2\beta \hbar ^{2} \\
2\beta \hbar ^{2} & 0%
\end{pmatrix}%
~,~V\left( x^{k},u^{C}\right) =V_{0}\Psi ^{2}-\beta \hbar ^{2}\Phi ^{2},
\label{pd.68}
\end{equation}%
therefore, the previous results apply. \ We show easily that the Ricci
Scalar of the $2-$dimensional metric $H_{AB}$ (\ref{pd.68}) vanishes, hence $%
H_{AB}$ is the two dimensional flat metric; furthermore, in this coordinate
system the connection coefficients are $C_{BC}^{A}=0$. \ Hence, the
Euler-Lagrange equations (\ref{pd.02}) for Lagrangian (\ref{pd.67}) are%
\begin{equation}
g^{\mu \nu }\Psi _{,\mu \nu }-\Gamma ^{\mu }\Psi _{,\mu }-\Phi =0
\label{pd.69}
\end{equation}%
\begin{equation}
2\beta \hbar ^{2}\left( g^{\mu \nu }\Phi _{,\mu \nu }-\Gamma ^{\mu }\Phi
_{,\mu }\right) +\left( V_{0}\Psi +\Phi \right) =0  \label{pd.70}
\end{equation}%
where equation (\ref{pd.69}) is the constraint\footnote{%
We remark that equations (\ref{pd.69}), (\ref{pd.70}) form a singular
perturbation system because $\beta \hbar ^{2}<<1$.} $\Phi =\Delta _{g}\Psi $.

$H_{AB}$ being a two dimensional flat space, admits a six dimensional Affine
algebra. In this coordinate system, the two gradient KVs of $H_{AB}~$are $%
K^{1}=\partial _{\Psi }$, $K^{2}=\partial _{\Phi },$ the non gradient KV is $%
R=2\beta \hbar ^{2}\Psi \partial _{\Psi }+\left( \Psi +2\beta \hbar ^{2}\Phi
\right) \partial _{\Phi }$; the gradient HV is $Y^{A}=\Psi \partial _{\Psi
}+\Phi \partial _{\Phi },$ and the proper ACs are$~A^{1}=\Psi \partial
_{\Psi },~A^{2}=\Phi \partial _{\Phi }~,~A^{3}=\Phi \partial _{\Psi }$ and $%
A^{4}=\Psi \partial _{\Phi }$. By replacing (\ref{pd.68}), in the constraint
equations of Theorems \ref{Lsym} and \ref{Nsym}, we get the following result.

\begin{corollary}
\label{CorGUP}The dynamical system with Lagrangian (\ref{pd.67}), which
describes the modified Klein-Gordon equation of a particle in GUP, admits as
Lie point symmetries the Killing vectors of the space of the independent
variables $g_{\mu \nu }$, plus the three vector fields $K_{G}=$ $b^{1}\left(
x^{k}\right) K^{1}+b^{2}\left( x^{k}\right) K^{2},~Y^{A}$ and $%
Z^{A}=A^{1}+2\beta \hbar ^{2}A^{3}-V_{0}A^{4}~$where $b^{1},~b^{2}$ solve
the system (\ref{pd.69}), (\ref{pd.70}). \ As far as the Noether point
symmetries of Lagrangian (\ref{pd.68}) are concerned, all the Lie point
symmetries of the dynamical system except the $Z^{A}$ are also Noether point
symmetries.
\end{corollary}

Furthermore, concerning the maximal dimension of the Lie algebra of the
system (\ref{pd.69})-(\ref{pd.70}), from corollary \ref{CorGUP} follows:

\begin{corollary}
\label{CorGUP2} If the system (\ref{pd.69})-(\ref{pd.70}) of $n$ independent
variables is invariant under a group of one parameter point transformations $%
G_{P}$, then $3\leq \dim G_{P}\leq \frac{n\left( n+1\right) }{2}+3;$ the
right equality holds if and only if the space of the independent variables
is a maximally symmetric space, and the left equality holds if and only if
the space of the independent variables does not admit a Killing vector field.
\end{corollary}

In the following, we apply theorem \ref{CorGUP} in order to determine the
Lie and the Noether point symmetries of the system (\ref{pd.69}), (\ref%
{pd.70}), in two cases of special interest: A) the four dimensional
Minkowski spacetime $M^{4},$ and B) the two dimensional hyperbolic sphere.

\subsection{Case A: The Minkowski spacetime $M^{4}$}

\label{M4ex}

In this case $g_{\mu \nu }=\eta _{\mu \nu }$ and Lagrangian (\ref{pd.67})
becomes%
\begin{eqnarray}
L_{M^{4}} &=&\frac{1}{2}\left( \left( \Psi _{,t}\right) ^{2}-\left( \Psi
_{,x}\right) ^{2}-\left( \Psi _{,y}\right) ^{2}-\left( \Psi _{z}\right)
^{2}\right) -\left( \frac{1}{2}V_{0}\Psi ^{2}-\beta \hbar ^{2}\Phi
^{2}\right) +  \notag \\
&&+2\beta \hbar ^{2}\left( \Psi _{,t}\Phi _{,t}-\Psi _{,x}\Phi _{,x}-\Psi
_{,y}\Phi _{,y}-\Psi _{,z}\Phi _{,z}\right) .  \label{pd.71}
\end{eqnarray}%
Equations (\ref{pd.69})-(\ref{pd.70}) are%
\begin{equation}
\Psi _{,tt}-\Psi _{,xx}-\Psi _{,yy}-\Psi _{,zz}-\Phi =0  \label{pd.72}
\end{equation}%
\begin{equation}
2\beta \hbar ^{2}\left( \Phi _{,tt}-\Phi _{,xx}-\Phi _{,yy}-\Phi
_{,zz}\right) +\left( V_{0}\Psi +\Phi \right) =0.  \label{pd.73}
\end{equation}

The $M^{4}$ spacetime admits a ten dimensional Killing algebra, hence from
theorem \ref{CorGUP}, we have that the system (\ref{pd.72})-(\ref{pd.73})
admits 13 Lie point symmetries, and the Lagrangian (\ref{pd.71}) admits 12
Noether point symmetries.

Using the zero order invariants of the Lie point symmetry vectors which span
the Lie algebra $\left\{ \partial _{y},\partial _{,z},\partial
_{t}+cY^{A}\right\} ,$ we find the invariant solutions for the wave function
$\Psi \left( t,x,y,z\right) $ of the system (\ref{pd.72})-(\ref{pd.73}) to
be
\begin{equation}
\Psi \left( t,x,y,z\right) =e^{ct}\left( c_{1}e^{\mu \left( \beta \hbar
^{2}\right) x}+c_{2}e^{-\mu \left( \beta \hbar ^{2}\right) x}+c_{3}e^{\nu
\left( \beta \hbar ^{2}\right) x}+c_{4}e^{-\nu \left( \beta \hbar
^{2}\right) x}\right)  \label{pd.73b}
\end{equation}%
where $\mu \left( \beta \hbar ^{2}\right) =\frac{1}{2\beta \hbar ^{2}}\sqrt{%
4c^{2}\left( \beta \hbar ^{2}\right) ^{2}+\beta \hbar ^{2}\left( 1-\lambda
\right) }$, $\nu \left( \beta \hbar ^{2}\right) =\frac{1}{2\beta \hbar ^{2}}%
\sqrt{4c^{2}\left( \beta \hbar ^{2}\right) ^{2}+\beta \hbar ^{2}\left(
1+\lambda \right) }~$\ and $\lambda =\sqrt{1-8V_{0}\beta \hbar ^{2}}.$
However, symmetries can also be used in order to transform solutions into
solutions. The Lie point symmetry $Z^{A}$ gives us the point transformation%
\begin{eqnarray}
\bar{\Psi}\left( \Psi ,\Phi ,\varepsilon \right) &=&\frac{\exp \left( \frac{%
\left( 1-\lambda \right) }{2}\varepsilon \right) }{2\lambda }\left[ \left(
4\Phi \beta \hbar ^{2}+\left( 1+\lambda \right) \Psi \right) e^{\varepsilon
\lambda }-\left( 4\beta \hbar ^{2}\Phi +\left( 1-\lambda \right) \Psi
\right) \right]  \label{pd.73c} \\
\bar{\Phi}\left( \Psi ,\Phi ,\varepsilon \right) &=&\frac{\exp \left( \frac{%
\left( 1-\lambda \right) }{2}\varepsilon \right) }{2\lambda }\left[ \left(
\left( \lambda -1\right) \Phi -2V_{0}\Psi \right) e^{\varepsilon \lambda
}+\left( \left( 1+\lambda \right) \Phi +2V_{0}\Psi \right) \right]
\label{pd.73d}
\end{eqnarray}%
that is, solution (\ref{pd.73b}) is transformed to the following solution
\begin{equation}
\bar{\Psi}\left( t,x,y,z\right) =e^{ct}\left( c_{1}^{\prime }e^{\mu \left(
\beta \hbar ^{2}\right) x}+c_{2}^{\prime }e^{-\mu \left( \beta \hbar
^{2}\right) x}+c_{3}^{\prime }e^{\nu \left( \beta \hbar ^{2}\right)
x}+c_{4}^{\prime }e^{-\nu \left( \beta \hbar ^{2}\right) x}\right)
\label{pd.73e}
\end{equation}%
where now the new constants $c_{\left( 1,2\right) }^{\prime }$ are%
\begin{equation}
c_{\left( 1,2\right) }^{\prime }=\exp \left( -\frac{1+\lambda }{2}%
\varepsilon \right) c_{\left( 1,2\right) }~,~c_{\left( 3,4\right) }^{\prime
}=\exp \left( -\frac{1-\lambda }{2}\varepsilon \right) c_{\left( 3,4\right)
}.  \label{pd.73f}
\end{equation}

\subsection{Case B: The Hyperbolic sphere $S^{2}$}

Consider now the two dimensional hyperbolic sphere with line element
\begin{equation}
ds^{2}=d\theta ^{2}-e^{2\theta }d\phi ^{2}.  \label{pd.74}
\end{equation}

In this space, Lagrangian (\ref{pd.67}) becomes%
\begin{eqnarray}
L_{\mathcal{S}} &=&\frac{e^{\theta }}{2}\left( \left( \Psi _{,\theta
}\right) ^{2}-e^{-2\theta }\left( \Psi _{\phi }\right) ^{2}\right) +  \notag
\\
&&+2\beta \hbar ^{2}e^{\theta }\left( \Psi _{,\theta }\Phi _{,\theta
}-e^{-2\theta }\Psi _{,\phi }\Phi _{,\phi }\right) +  \notag \\
&&-e^{\theta }\left( \frac{1}{2}V_{0}\Psi ^{2}-\beta \hbar ^{2}\Phi
^{2}\right) .  \label{pd.75a}
\end{eqnarray}

Hence in the two dimensional hyperbolic sphere (\ref{pd.74}), the wave
function of the particle is computed from the following system of equations
\begin{equation}
\Psi _{,\theta \theta }-e^{-2\theta }\Psi _{,\phi \phi }+\Psi _{,\theta
}-\Phi =0  \label{pd.76}
\end{equation}%
\begin{equation}
2\beta \hbar ^{2}\left( \Phi _{,\theta \theta }-e^{-2\theta }\Phi _{,\phi
\phi }+\Phi _{,\theta }\right) +\left( V_{0}\Psi +\Phi \right) =0.
\label{pd.77}
\end{equation}

The two dimensional sphere (\ref{pd.74}) admits a three dimensional Killing
algebra, the $SO\left( 3\right) $. Hence, from Corollary \ref{CorGUP} we
have that the system (\ref{pd.76})-(\ref{pd.77}) admits 6 Lie point
symmetries, and the Lagrangian (\ref{pd.75a}) admits 5 Noether point
symmetries.

The elements of $SO\left( 3\right) $ algebra in the coordinates of (\ref%
{pd.74}) are%
\begin{equation}
X^{1}=\partial _{\phi }~,~X^{2}=\partial _{\theta }-\phi \partial _{\phi }
\label{pd.78}
\end{equation}%
\begin{equation}
X^{3}=2\phi \partial _{\theta }-\left( \phi ^{2}+e^{-2\theta }\right)
\partial _{\phi }.  \label{pd.79}
\end{equation}

From the application of the Lie point symmetry $X^{1}+\alpha Y^{A}$ we have
the solution%
\begin{eqnarray}
\Psi _{1}\left( \theta ,\phi \right) &=&e^{\alpha \phi }e^{-\frac{\theta }{2}%
}\left[ b_{1}K_{\bar{\mu}}\left( \alpha e^{-\theta }\right) +b_{2}I_{\bar{\mu%
}}\left( \alpha e^{-\theta }\right) \right]  \notag \\
&&+e^{\alpha \phi }e^{-\frac{\theta }{2}}\left[ b_{3}K_{\bar{\nu}}\left(
\alpha e^{-\theta }\right) +b_{4}I_{\bar{\nu}}\left( \alpha e^{-\theta
}\right) \right]  \label{pd.80}
\end{eqnarray}%
where~$\bar{\mu}=$ $-\frac{\sqrt{\left( \beta \hbar ^{2}-1-\lambda \right) }%
}{2\sqrt{\beta \hbar ^{2}}}$,~$\bar{\nu}=-\frac{\sqrt{\left( \beta \hbar
^{2}-1+\lambda \right) }}{2\sqrt{\beta \hbar ^{2}}}$, $\lambda =\sqrt{%
1-8V_{0}\beta \hbar ^{2}},$ and $I_{l},K_{l}$ are the modified Bessel
functions of the first and second kind respectively.

By using the Lie point symmetry $X^{2}+\kappa Y^{A}$ we find the solution%
\begin{eqnarray}
\Psi _{2}\left( \theta ,\phi \right) &=&\left( \phi ^{2}e^{2\theta
}+1\right) ^{-\frac{\kappa }{2}}e^{\kappa \theta }\left[ b_{1}P_{\left( -%
\bar{\mu}-\frac{1}{2}\right) }^{\kappa }\left( \phi e^{\theta }\right)
+b_{2}Q_{\left( -\bar{\mu}-\frac{1}{2}\right) }^{\kappa }\left( \phi
e^{\theta }\right) \right] +  \notag \\
&&+\left( \phi ^{2}e^{2\theta }+1\right) ^{-\frac{\kappa }{2}}\left[
b_{3}P_{\left( -\bar{\nu}-\frac{1}{2}\right) }^{\kappa }\left( \phi
e^{\theta }\right) +b_{4}Q_{\left( -\bar{\nu}-\frac{1}{2}\right) }^{\kappa
}\left( \phi e^{\theta }\right) \right]  \label{pd.81}
\end{eqnarray}%
where $P_{l}^{k},~Q_{l}^{k}$ are the associated Legendre functions of the
first and second kind respectively.

Furthermore, from the Lie symmetry $X^{3}+\sigma Y^{A}$ we find the solution%
\begin{eqnarray}
\Psi _{3}\left( \theta ,\phi \right) &=&\exp \left( \frac{\sigma \phi
e^{2\theta }}{\phi ^{2}e^{2\theta }+1}\right) \left[ b_{1}K_{\bar{\mu}%
}\left( \frac{\sigma e^{2\theta }}{\phi ^{2}e^{2\theta }+1}\right) +b_{2}I_{%
\bar{\mu}}\left( \frac{\sigma e^{2\theta }}{\phi ^{2}e^{2\theta }+1}\right) %
\right] +  \notag \\
&&+\exp \left( \frac{\sigma \phi e^{2\theta }}{\phi ^{2}e^{2\theta }+1}%
\right) \left[ b_{3}K_{\bar{\nu}}\left( \frac{\sigma e^{2\theta }}{\phi
^{2}e^{2\theta }+1}\right) +b_{4}I_{\bar{\nu}}\left( \frac{\sigma e^{2\theta
}}{\phi ^{2}e^{2\theta }+1}\right) \right]  \label{pd.82}
\end{eqnarray}

Finally, if we apply the point transformation (\ref{pd.73c})-(\ref{pd.73d}),
which follows from the Lie point symmetry $Z^{A},~$to the solutions (\ref%
{pd.80})-(\ref{pd.82}), we find that the new solutions are again (\ref{pd.80}%
)-(\ref{pd.82}), where now the constants $b_{1-4}\rightarrow $ $\bar{b}%
_{1-4} $ \ are given by the expressions (\ref{pd.73f}).

\section{Conclusion}

\label{conclusions}

The quasilinear systems of second-order differential equations ($n\times m$
systems) describe many important equations of relativistic Physics, both at
the classical as well as at the quantum level. Therefore, the study of Lie
and the Noether point symmetries of these equations is important in order to
establish invariant solutions and find conservation laws.

In this work we followed the geometric approach we have applied in our
previous studies for the point symmetries of $1\times m$ and $n\times 1$
systems in order to generalize it to the case of $n\times m$ systems. We
have shown that for the $n\times m$ systems of the form of (\ref{pd.02}),
with $n\geq 2$, the point symmetries follow from the CKVs and the ACs of the
underlying geometries of the $n~$independent and $\ $the $m~$dependent
variables. This result is consistent with the results concerning the $%
1\times m$ and $n\times 1$ systems. Moreover, for the Lagrangian (\ref{pd.01}%
) we derived the general form of the Noether symmetry vector and of the
Noetherian conservation laws. Specifically, we proved that the Noether point
symmetries are generated by the CKVs and the HVs of the metrics $g_{ij}$ and
$H_{AB}$ respectively.

We applied the above general results to a system of quasilinear Laplace
equations (which contains the geodesic equations and the wave equation as
special cases) and determined the Lie and the Noether point symmetries in
flat space and in spaces of constant non-vanishing curvature. By using
results of Differential Geometry we found the maximum dimension of the Lie
algebra of a system of quasilinear Laplace equations. In particular, we
showed that if the $n\times m$ system of quasilinear Laplace equations, with
$n>2$, admits $\frac{\left( n+2\right) \left( n+1\right) }{2}+m\left(
m+1\right) $ linear independent Lie point symmetries, then there exists a
coordinate system in which this system takes the simplest form $\delta
^{ij}u_{,ij}^{A}=0$.

A further application concerns the Klein-Gordon equation in GUP which is a
fourth-order partial differential equation. By using a Lagrange multiplier
we reduced this fourth-order equation to a system of two coupled
second-order equations in the form of the system (\ref{pd.02}), in which the
main Theorems of this work applies. \ We showed that the minimum Lie algebra
of that system is of dimension three. Furthermore, we determined the Lie and
the Noether point symmetries in a Minkowski space and in the 2-dimensional
hyperbolic space and subsequently we used the former in order to determine
invariant solutions of the wave function of a spin-$0$ particle.

{\large {\textbf{Acknowledgements}}} \newline
The research of AP was supported by FONDECYT postdoctoral grant no. 3160121.

\end{document}